\documentclass[english]{article}
\usepackage[T1]{fontenc}
\usepackage[latin9]{inputenc}
\usepackage{amsmath}
\usepackage{amssymb}
\usepackage{graphicx}

\makeatletter

\providecommand{\tabularnewline}{\\}

\usepackage[english]{babel}

\makeatother

\usepackage{babel}
\begin{document}

\title{Monte Carlo approximations of the Neumann problem }

\author{Sylvain Maire%
\thanks{Laboratoire des Sciences de l'Information et des Systemes (LSIS),
UMR6168, ISITV, Université de Toulon et du Var, Avenue G. Pompidou,
BP 56, 83262 La Valette du Var cedex, France \texttt{sylvain.maire@univ-tln.fr}%
} $\,\:$and Etienne Tanré%
\thanks{INRIA, EPI Tosca, 2004 route des Lucioles BP93 F-06902 Sophia-Antipolis
- France\protect \\
\texttt{Etienne.Tanre@inria.fr}%
}}
\date{August 21, 2013}
\maketitle
\begin{abstract}
We introduce Monte Carlo methods to compute the solution of elliptic
equations with pure Neumann boundary conditions. We first prove that
the solution obtained by the stochastic representation has a zero
mean value with respect to the invariant measure of the stochastic
process associated to the equation. Pointwise approximations are computed
by means of standard and new simulation schemes especially devised
for local time approximation on the boundary of the domain. Global
approximations are computed thanks to a stochastic spectral formulation
taking into account the property of zero mean value of the solution.
This stochastic formulation is asymptotically perfect in terms of
conditioning. Numerical examples are given on the Laplace operator
on a square domain with both pure Neumann and mixed Dirichlet-Neumann
boundary conditions. A more general convection-diffusion equation
is also numerically studied.
\end{abstract}
\textbf{Keywords.} Neumann Problem, Monte Carlo Methods, Spectral Methods,
Local Time Approximation
\section{Introduction}

The aim of this paper is to compute pointwise and global numerical
approximations of the solutions of pure Neumann problems for elliptic
equations by means of Monte Carlo procedures. The pointwise solutions
will be obtained via Feynman-Kac representations and the global ones
by stochastic spectral formulations inspired from the methods developed
in our previous works \cite{maire-tanre-2008,maire-tanre-2009}.

We first consider the Neumann problem for the Laplace operator in
a domain $D\subset\mathbb{R}^{d}$ which writes 
\begin{align*}
\begin{cases}
\frac{1}{2}\Delta u(x) & =-f(x)\quad\text{for}\quad x=(x_{1},\ldots,x_{d})\in D\\
\frac{1}{2}\frac{\partial u}{\partial n}(y) & =g(y)\quad\text{for}\quad y\in\partial D,
\end{cases}
\end{align*}
where $\frac{\partial u}{\partial n}$ stands for the incoming normal
derivative of $u$ on the boundary of the domain $D.$ Because the
boundary conditions are of Neumann type everywhere on the boundary,
this equation has a solution up to an additive constant whenever an
additional compatibility condition on $f$ and $g$ is satisfied (see
\cite{giraud1933}). This compatibility condition writes 
\[
\int_{D}f(x)dx+\int_{\partial D}g(y)dy=0.
\]
In the case of the Poisson equation with Dirichlet boundary conditions,
the probabilistic representations of the solutions are based on Feynman-Kac
formulae which are stopped at the first hitting time of the boundary
by a Brownian motion. In the case of mixed boundary conditions, that
is of Neumann type on a part of the boundary and of Dirichlet type
on the other part, the (reflected) Brownian motion is still stopped
at the first hitting time of the Dirichlet boundary. In the case of
pure Neumann boundary conditions the Brownian motion never stops and
the probabilistic representation introduced by Brosamler \cite{brosamler}
in the case $f=0$ is the limit in time of the Feynman-Kac representation
of the solution of the relative Cauchy problem. We have 
\[
u(x)=\lim_{T\rightarrow\infty}\mathbb{E}_{x}\left(\int_{0}^{T}g(B_{s})dV_{s}\right)
\]
 where $(B_{t},t\geq0)$ stands for the reflected Brownian motion
in $D$ and $(V_{t},t\geq0)$ is its associated local time on the
boundary. The proof of this representation is based on probabilistic
potential theory combined with a representation of additive functionals.
This representation has been extended by Benchérif-Madani and Pardoux
\cite{bencherif_pardoux} to more general advection-diffusion equations
\[
Lu(x)=-f(x),
\]
 with boundary conditions $\frac{\partial u}{\partial n_{a}}(y)=g(y)$,
where 
\[
Lu(x)=\frac{1}{2}\sum_{i,j=1}^{d}a_{i,j}(x)\frac{\partial^{2}u(x)}{\partial x_{i}\partial x_{j}}+\sum_{i=1}^{d}b_{i}(x)\frac{\partial u(x)}{\partial x_{i}}
\]
 and $\frac{\partial u}{\partial n_{a}}(.)=\frac{1}{2}a\nabla(.).n$
stands for the conormal derivative. Some regularity assumptions are
required on the domain $D,$ the functions $f$ and $g$, the coefficients
$a_{i,j}$ and $b_{i}$ and also uniform ellipticity conditions and
boundness for the symmetric matrix $a.$ We have 
\begin{equation}
u(x)=\lim_{T\rightarrow\infty}u(x,T)=\lim_{T\rightarrow\infty}\mathbb{E}{}_{x}\left(\int_{0}^{T}f(X_{s})ds+\int_{0}^{T}g(X_{s})dV_{s}\right)\label{eq:representation}
\end{equation}
 where $(X_{t},t\geq0)$ is the reflected diffusion process and $u(x,T)$
the solution of the Cauchy problem associated to the operator $L$.
The compatibility condition becomes 
\[
\int_{D}f(x)p(x)dx+\int_{\partial D}g(y)p(y)dy=0
\]
 where $p(x)$ is the density of the invariant measure associated
to $(X_{t},t\geq0).$ The proof of this more general representation
is now based on exponential ergodicity.

We intend to use the previous stochastic representation for the numerical
computation of the Neumann problem by means of Monte Carlo procedures.
We need to overcome several problems. To compute the pointwise solution,
we have to choose a time $T_{0}$ to stop the trajectories, deal efficiently
with the approximation of the Neumann boundary conditions (which involves
local time approximation) and understand which additive constant is
obtained using this representation. The Monte Carlo approximation
of the Neumann problem associated to the Laplace operator has already
been treated using the walk on the boundary method \cite{sabsimonov}.
The truncation problems of the divergent Neumann serie linked to this
walk are very similar to the choice of the stopping time $T_{0}$
of our trajectories. This method deals very naturally with the Neumann
boundary conditions but may suffer of an increase of the variance
for the computation of the solution at points close to the boundary
of the domain. Some estimators with a reduced variance are also proposed
in \cite{sabsimonov}. They require to find the part of the boundary
where the boundary term is positive. To obtain the global solution,
we need to choose and to build carefully the basis functions involved
in the stochastic spectral formulations in order to obtain an unique
and well characterized solution. The rest of the paper is organized
as follows.

In section \ref{sec:properties_stochastic_representation}, we first
prove that the representation of Benchérif-Madani and Pardoux \cite{bencherif_pardoux}
naturally leads to a solution $u(x)$ such that its mean value with
respect to the invariant measure is zero. For a practical computation,
we need to replace the limit in representation (\ref{eq:representation})
by a finite time $T_{0}.$ This introduces a bias that will be linked
to the second eigenvalue of the operator $L$ with pure Neumann boundary
conditions. The question of the variance of our estimators as a function
of $T_{0}$ is also studied. We remark that the variance has a linear
growth as a function of $T_{0}.$

In section \ref{sec:Euler-scheme}, we describe a general algorithm
to approximate the solution of general elliptic Neumann problems by
means of the Euler scheme coupled with a local time approximation
method developed in \cite{jiang-knight}.

Section \ref{sec:WOS} is devoted to walk on spheres algorithms which
can be used in the special case of the Poisson equation. We describe
how to adapt them to a Cauchy problem with a finite horizon $T_{0}.$
Furthermore, two new schemes with an increased order of convergence
are introduced to deal with inhomogeneous Neumann boundary conditions.

In section \ref{sec:Numerical-results}, we give some numerical results
on mixed Dirichlet-Neumann problems and on pure Neumann problems with
different degrees of difficulty. For the pure Neumann problem, the
scheme that is used to compute the representation modifies the additive
constant from which the solution depends.

Finally in section \ref{sec:Stochastic-spectral-methods}, we develop
stochastic spectral methods in order to obtain a very accurate and
global approximation of the solution. In addition to the general methodology
developed in \cite{maire-tanre-2008,maire-tanre-2009}, we need to
develop centering procedures for our approximation basis in order
to project our solution in the space with zero mean value with respect
to $p(x).$ Numerical examples show that our method is very efficient
in terms of accuracy and conditioning in the both cases where the
invariant measure is known analytically or only approximated.

\section{Some properties of the stochastic representation\label{sec:properties_stochastic_representation}}

For theoretical purpose, representation (\ref{eq:representation})
is useful. For a numerical purpose, we need to understand which solution
is effectively computed (not only up to an additive constant) and
the consequences of replacing the limit by a fixed time $T_{0}$.
We shall prove in the following that the solution $u(x,T)$ has asymptotically
a mean value zero with respect to the invariant measure and that its
variance is increasing mainly linearly as a function of $T.$

We consider a general diffusion operator $L$ verifying the properties
mentioned in the introduction \cite{bencherif_pardoux}. We denote
by $(X_{t},t\geq0)$ its associated reflected diffusion process

\begin{equation}
X_{t}=x+\int_{0}^{t}b(X_{s})ds+\int_{0}^{t}\sigma(X_{s})dW_{s}-\frac{1}{2}\int_{0}^{t}\gamma(X_{s})dV_{s}\label{eq:dynamiquex}
\end{equation}
 where the matrix $\sigma$ is such that $\sigma\sigma^{T}=a$, $(W_{t},t\geq0)$
is a $d$ dimensional Brownian motion, $(V_{t},t\geq0)$ is the local
time on the boundary and $\gamma$ is in the conormal direction 
\[
\gamma(x)=an(x).
\]
 We then define the solution $u(x,T)$ of the Cauchy problem 
\[
u(x,T)=\mathbb{E}{}_{x}\left(\int_{0}^{T}f(X_{s})ds+\int_{0}^{T}g(X_{s})dV_{s}\right)
\]
 where $f$ and $g$ are bounded and verify the compatibility condition
\[
\int_{D}f(x)p(x)dx+\int_{\partial D}g(y)p(y)dy=0.
\]
 Following the arguments developed in \cite{bencherif_pardoux}, there
exists a bounded function $f_{1}$ such that 
\[
u(x,T)=\mathbb{E}_{x}\left(\int_{0}^{T}f_{1}(X_{s})ds\right)
\]
 with the new compatibility condition 
\[
\int_{D}f_{1}(x)p(x)dx=0.
\]

\subsection{Characterization of the solution}

We can first notice that 
\[
\int_{D}\mathbb{E}_{x}(f_{1}(X_{s}))p(x)dx=\mathbb{E}_{\mu}(f_{1}(X_{s}))=\int_{D}f_{1}(x)p(x)dx=0
\]
 for any time $s$ because $\mu$ is the invariant law for \eqref{eq:dynamiquex}.
Using the Fubini theorem, this leads to 
\[
\int_{D}u(x,T)p(x)dx=0
\]
 for any time $T$. In the case of the Laplace operator, it has been
proven in \cite{bass-hsu} that the convergence in time of $u(x,T)$
towards $u(x)$ is uniform in $x$ and thus the desired property holds.
This uniform convergence relies on special properties of the Brownian
motion that cannot be extended so easily to a general diffusion process.
These properties are a bound on the transition density function $p(t,x,y)$
of the reflecting Brownian motion and its spectral expansion \cite{burdzy}.
This spectral expansion writes 
\[
p(t,x,y)=\frac{1}{\left|D\right|}+\sum_{n=1}^{\infty}\exp(-\lambda_{n}t)\phi_{n}(x)\phi_{n}(y)
\]
 where $\frac{1}{\left|D\right|}$ is the uniform measure in $D$,
$(\lambda_{n},n\geq0)$ and $(\phi_{n},n\geq0)$ the eigenelements
of $-\frac{1}{2}\Delta$. The first term $\frac{1}{\left|D\right|}$
corresponds to the first eigenvalue $\lambda_{0}=0$ and 
\[\phi_{0}(x)=\frac{1}{\sqrt{D}}.\]
The remaining eigenvalues $(\lambda_{n},n\geq1)$ are strictly positive.
Theorem 2.4 of \cite{bass-hsu} says that $p(t,x,y)$ converges exponentially
fast and uniformly to $\frac{1}{\left|D\right|}$ at a speed linked
to the second eigenvalue $-\lambda_{1}$ of $\frac{1}{2}\Delta.$
This can guide us for the choice of the time $T_{0}$ at which we
should stop our trajectories.

In the general case, we have for any $T\geq0,$ 
\[
\left|\int_{D}u(x)p(x)dx\right|^{2}=\left|\int_{D}u(x)p(x)dx-\int_{D}u(x,T)p(x)dx\right|^{2}\leq\int_{D}(u(x)-u(x,T))^{2}p(x)dx
\]
 thanks to the Cauchy-Schwarz inequality. In lemma 3 of \cite{bencherif_pardoux},
it is proven an exponential ergodicity for the solution $w(x,t)=\mathbb{E}_{x}(f_{1}(X_{t}))$
of the parabolic problem with an initial condition $f_{1}$ satisfying
the centering condition, no source term and homogeneous Neumann boundary
conditions. This is summarized in the inequality 
\[
\int_{D}w(x,t)^{2}p(x)dx\leq\exp(-rt)\int_{D}f_{1}(x)^{2}p(x)dx
\]
 where $r$ is a positive constant. We can write 
\[
\int_{D}(u(x)-u(x,T))^{2}p(x)dx\leq\int_{D}p(x)dx\int_{T}^{\infty}\exp(\frac{r}{2}s)w(x,s)^{2}ds\times\int_{T}^{\infty}\exp(-\frac{r}{2}s)ds
\]
 thanks again to the Cauchy-Schwarz inequality and finally 
\[
\left|\int_{D}u(x)p(x)dx\right|\leq\frac{2}{r}\exp(-rT)\int_{D}f_{1}(x)^{2}p(x)dx
\]
 thanks to the Fubini theorem and to the exponential ergodicity bound.
This proves that $\int_{D}u(x)p(x)dx=0.$

\subsection{Estimation of the variance}

In this section, we estimate the variance of $\int_{0}^{T}f_{1}(\bar{B}_{s})ds$
as a function of $T$ in the simplified case of the Poisson equation
in dimension one in the interval $D=[a,b]$ with boundary conditions
$u^{'}(a)=u^{'}(b)=0$. We have 
\[
Var\left(\int_{0}^{T}f_{1}(\bar{B}_{s})ds\right)=\mathbb{E}_{x}\left(\int_{0}^{T}f_{1}(\bar{B}_{s})ds\right)^{2}-\left(\mathbb{E}_{x}\int_{0}^{T}f_{1}(\bar{B}_{s})ds\right)^{2}
\]
 where $(\bar{B}_{s},s\geq0)$ is a reflected Brownian motion in $D$.
As 
\[
\lim_{T\rightarrow\infty}\left(\mathbb{E}_{x}\int_{0}^{T}f_{1}(\bar{B}_{s})ds\right)^{2}=u^{2}(x),
\]
 the second term in the variance is obviously bounded because $u$
is a continuous function in a bounded domain. For the first term,
thanks to Itô-Tanaka formula we have 
\[
u(\bar{B}_{T})=u(x)+\int_{0}^{T}Lu(\bar{B}_{s})ds+\int_{0}^{T}u^{'}(\bar{B}_{s})dB_{s}=u(x)-\int_{0}^{T}f_{1}(\bar{B}_{s})ds+\int_{0}^{T}u^{'}(\bar{B}_{s})dB_{s}.
\]
 Thanks to the boundary conditions, the usual term involving the local
time process is equal to zero. The above formula leads to 
\begin{eqnarray*}
\left(\int_{0}^{T}f_{1}(\bar{B}_{s})ds\right)^{2} & = & \left(u(x)-u(\bar{B}_{T})\right)^{2}+\left(\int_{0}^{T}u^{'}(\bar{B}_{s})dB_{s}\right)^{2}\\
 &  & +2\left(u(x)-u(\bar{B}_{T})\right)\left(\int_{0}^{T}u^{'}(\bar{B}_{s})dB_{s}\right).
\end{eqnarray*}
 Using that the invariant measure is uniform in $D$, we have 
\[
\lim_{T\rightarrow\infty}\mathbb{E}_{x}\left[\left(u(x)-u(\bar{B}_{T})\right)\right]^{2}=\frac{1}{\left|D\right|}\int_{D}\left(u(x)-u(y)\right)^{2}dy
\]
 which proves that this term is bounded and its limit identified.
Then from the isometry property 
\[
\mathbb{E}_{x}\left[\left(\int_{0}^{T}u^{'}(\bar{B}_{s})dB_{s}\right)^{2}\right]=\mathbb{E}_{x}\left[\int_{0}^{T}\left(u^{'}(\bar{B}_{s})\right)^{2}ds\right]
\]
 and from the ergodic theorem 
\[
\mathbb{E}_{x}\left[\left(\int_{0}^{T}u^{'}(\bar{B}_{s})dB_{s}\right)^{2}\right]=T\mathbb{E}_{x}\left[\frac{1}{T}\int_{0}^{T}\left(u^{'}(\bar{B}_{s})\right)^{2}ds\right]\simeq T\frac{1}{\left|D\right|}\int_{D}\left(u^{'}(y)\right)^{2}dy.
\]
 Thanks now to the Cauchy-Schwarz inequality and to the previous
bounds

\[
\left|\mathbb{E}_{x}\left(u(x)-u(\bar{B}_{T})\right)\left(\int_{0}^{T}u^{'}(\bar{B}_{s})dB_{s}\right)\right|\leq C\sqrt{T}\left(\frac{1}{\left|D\right|}\int_{D}\left(u^{'}(y)\right)^{2}dy\right)^{1/2}
\]
 which proves finally that asymptotically

\[
C_{1}-C_{2}\sqrt{T}+C_{3}T\leq Var(\int_{0}^{T}f(\bar{B}_{s})ds)\leq C_{1}+C_{2}\sqrt{T}+C_{3}T
\]
 for some constants $C_{1},C_{2}$ and where $C_{3}=\frac{1}{\left|D\right|}\int_{D}\left(u^{'}(y)\right)^{2}dy$.
The lower and upper bounds on the variance show that it increases
essentially as a linear function of $T$ (except if $u(x)$ is constant).
This proves that the time chosen to stop the trajectories is both
crucial in terms of bias and in terms of variance.

\section{Euler scheme approximations\label{sec:Euler-scheme}}

We now present a scheme used to approximate the solution $u$ inspired
by the representation (\ref{eq:representation}) but where the limit
is replaced by a fixed finite time $T_{0}$. This method relies on
weak approximations. First of all, in a general setting, we have to
approximate the reflected diffusion process $(X_{t},0\leq t\leq T_{0})$
associated to the infinitesimal generator $L$ in the domain $D\subset\mathbb{R}^{d}$.
The approximation of diffusion processes is usually done using the
Euler scheme. The errors and the error expansions on the weak approximations
using this scheme in the whole space have been first studied in \cite{talay-tubaro}.
The weak approximation of diffusions with absorbing (Dirichlet) and
reflected (Neumann) boundary conditions have been treated respectively
in \cite{gobet_2000} and \cite{bossy-gobet-talay-2004}. Given a
time step $\delta$, the approximation of $(X_{t},t\geq0)$ by a standard
reflected Euler scheme $(\bar{X}_{k\delta},0\leq k\leq T_{0}/\delta)$
can be described by the following procedure. 
\begin{enumerate}
\item for all $i=1..d,$ 
\[
\tilde{X}_{(k+1)\delta}^{i}=\bar{X}_{k\delta}^{i}+b_{i}(\bar{X}_{k\delta})\delta+\sum_{j}\sigma_{i,j}(\bar{X}_{k\delta})(W_{(k+1)\delta}^{j}-W_{k\delta}^{j}).
\]

\item if $\tilde{X}_{(k+1)\delta}\in D$, we set $\bar{X}_{(k+1)\delta}=\tilde{X}_{(k+1)\delta},$ 
\item else, if $\tilde{X}_{(k+1)\delta}\notin D$, we have to choose another
position inside the domain $D$. The final position $\bar{X}_{(k+1)\delta}$
is the symmetrized of the Euler scheme $\tilde{X}_{(k+1)\delta}$
in the conormal direction $a.n$ (see \cite{bossy-gobet-talay-2004}). 
\end{enumerate}
To compute our representation, we first need to approximate the integral
$\int_{0}^{T_{0}}f(X_{s})ds.$ This can be done using the rectangle
method by the discrete sum 
\[
\delta\sum_{k=0}^{\frac{T_{0}}{\delta}-1}f(\bar{X}_{k\delta}).
\]
The approximation of the term $\int_{0}^{T_{0}}g(X_{s})dV_{s}$ is
less classical (see e.g. \cite{jiang-knight}). We introduce
a delocalization parameter $\xi$ and compute the approximation 
\[
\delta\sum_{k=0}^{\frac{T_{0}}{\delta}-1}g(\pi(\bar{X}_{k\delta}))K_{\xi}(\bar{X}_{k\delta},\pi(\bar{X}_{k\delta}))
\]
 where $K_{\xi}$ stands for instance for the Gaussian kernel 
\[
K_{\xi}(x,y)=\frac{1}{(2\pi)^{d/2}\xi^{d}}\exp\left(-\sum_{i=1}^{d}\frac{\left(x_{i}-y_{i}\right)^{2}}{\xi^{2}}\right).
\]
 We have also used the orthogonal projection $\pi$ on the boundary
$\partial D$. The delocalization parameter must be chosen carefully.
We will observe in section \ref{sec:Numerical-results} that it should
not be too small unless the variance becomes very large. In this section,
we will also consider a problem with mixed boundary conditions. The
treatment of the Dirichlet boundary conditions will rely on half-space
approximations \cite{gobet-2001}. An additional stopping test of
the trajectories based on a Brownian bridge enables to obtain an order
1 weak error in $\delta$ on the part of the boundary with Dirichlet
boundary conditions. The layer method described in \cite{milstre}
can also be used. It is especially adapted to parabolic problems and
enables to obtain very accurate approximations of such problems. It
relies however on a Markov chain which is more difficult to simulate
than the Euler scheme.

\section{Walk on spheres approximations\label{sec:WOS}}

\subsection{Introduction}

The Euler scheme can be used for the simulation of a wide class of
stochastic processes linked to second order elliptic operators. If
we deal with the Laplace operator or with divergence form operators
with constant or piecewise constant diffusion coefficients \cite{lejay-maire},
more efficient simulations of Brownian paths are available like walk
on spheres (WOS) \cite{sabelfeld} or walk on rectangles algorithms
\cite{deaconu-lejay}.

In the case of Dirichlet boundary conditions, the stochastic representation
of the solution implies a Brownian path up to the first time it hits
the boundary $\partial D$ of the domain. To solve the Neumann problem,
we should use the same scheme as soon as the particle is away from
the boundary, that is it has not reached the $\varepsilon$-absorption
layer. Then, we replace the particle inside the domain and run again
the same scheme until hitting one more time the boundary and so on.
Some efficient ways to replace the particle inside the domain are
described in section \ref{sub:boundarycond}.

Let us now describe the WOS method to compute the solution of the
Laplace equation at a point $x\in D$. The walk starts at $x$ and
jumps from one sphere to another until it reaches an $\varepsilon$-absorption
layer ($=\left\{ y\in D,d(y,\partial D)\leq\varepsilon\right\} $). The
spheres are built so that the jumps are as large as possible by taking
the radius of the next sphere as the distance to the boundary $\partial D.$
The next point is the first hitting point of the sphere by a Brownian
motion started at the center. So, thanks to the isotropy of the Brownian
motion, it has to be chosen uniformly on this sphere. We stop the
walk at the first time the selected point on the sphere belongs to
the $\varepsilon$-absorption layer. The contribution of each walk to
the solution is the value of the boundary term at a projection (generally
the orthogonal one) of the current position on the boundary. Finally,
the approximate solution is the average of the contributions of walks.
In the case of the Poisson equation, it is also possible to compute
the contribution of the source term during this walk \cite{hwang-etal}
using for each sphere $S$ of radius $r_{n}$ and center $x_{n}$
a Green function $G$ conditioned by the exit point $z$ writing 
\begin{equation}
\mathbb{E}_{x_{n}}\left(\int_{0}^{\tau_{S}}f(B_{t})dt\,\big|\, B_{\tau_{S}}=z\right)=\frac{r_{n}^{2}}{d}\int_{S}G(z,y)f(y)dy.\label{eq:loicondsphere}
\end{equation}
 At least in dimension 2, the conditional density is known exactly
and it is possible to sample easily from it. Thus we can obtain the
contribution of the source term to the whole trajectory.

\subsection{WOS for the pure Neumann problem\label{sub:WOS-pure-neumann}}

Regardless the problem of making a good choice of a final time $T_{0}$
to stop our trajectories, our task is to deal efficiently with the
inhomogeneous Neumann boundary conditions while simulating accurately
the elapsed time and position from the start of the trajectories.
The description of the tools used when the walk hits the boundary
is given in the next subsection. We need to adapt the WOS method to
the case of a finite horizon $T_{0}$. First, formula (\ref{eq:loicondsphere})
gives not the exit time $\tau_{S}$ but only the average time $\frac{r_{n}^{2}}{d}$
spent in the sphere $S$ of radius $r_{n}.$ Second, we want to approximate
$\int_{0}^{\tau_{S}}f(B_{t})dt$ knowing that $B_{\tau_{S}}=z.$ It
is possible to sample from $\tau_{S}$ using for instance the inverse
method as its law has a closed form given by a serie. However, we
need here to compute more complicated quantities involving positions
of the Brownian motion before it hits the boundary conditioned by
the value of $\tau_{S}.$ We actually do not know closed forms for
these quantities so we sample from them using precomputed discretized
trajectories as described as follows. Note that anyway sampling from
a law that is given by a serie by the inverse method requires an iterative
method which may be time consuming compared to precomputation. The
previous integral can be transformed writing 
\[
\int_{0}^{\tau_{S}}f(B_{t})dt=\tau_{S}\mathbb{E}^{U}\left[f\left(B_{U\tau_{S}}\right)\right]
\]
where $U$ is an uniform r.v. on $[0,1]$ independent of $\left(B_{t},0\leq t\leq\tau_{S}\right)$
and $\mathbb{E}^{U}$ stands for the expectation with respect to $U$.
This transformation called the one random point method has been used
in our papers \cite{maire-tanre-2008,maire-tanre-2009} and also in
\cite{reutenauer-tanre-2008} in the context of financial mathematics.
Its interest is to give a lot cheaper evaluation of the integral for
only a small increase of the variance. 

Thanks to the isotropy of the Brownian Motion, we only need to simulate
the couple $(\tau_{S_{1}},B_{U\tau_{S_{1}}})$ knowing that $B_{\tau_{S_{1}}}=(1,0)$
in the unit circle. To do this, we simulate an absorbed Brownian motion
starting at the center of the unit circle using the Euler scheme with
a very small time step $\delta$ and the half space approximation
\cite{gobet-2001}. The walk stops after a random number $M$ of steps
at an exit point $z=\exp(i\theta).$ The exit time $\tau_{S}$ is
approximated by $M\delta$, a point $(x,y)$ is picked uniformly at
random among the $M$ points of the discretized trajectory and is
rotated using a rotation of angle $-\theta.$ We can thus obtain the
empirical joint law of $(\tau_{S_{1}},B_{U\tau_{S_{1}}})$ knowing
that $B_{\tau_{S_{1}}}=(1,0).$ Samples from this empirical law are
precomputed and stored in a large list. To sample from the couple
$(\tau_{S_{1}},B_{U\tau_{S_{1}}})$, we pick uniformly at random a
couple (exit time, position) in the list and an angle $\alpha$ uniformly
in $[0,2\pi]$ to rotate the position of the previous couple. 

Dealing with the pure Neumann problem, we have to take care to another
difference with the resolution of the Dirichlet problem: at the last
step, the final time $T_{0}$ is before the exit time of the sphere.
Assuming that time is reinitialized to zero before the last sphere,
we have to compute integrals of the form $\int_{0}^{T_{1}}f(B_{s})ds$
knowing that $B_{\tau_{S}}=z$ and $T_{1}\leq\tau_{S}.$ So, we also
need to store $Q$ full trajectories (time, position) before the exit
time $\tau_{S}$ to approximate the above integrals.

\subsection{Replacement after hitting the boundary\label{sub:boundarycond}}

It remains to deal with the inhomogeneous Neumann boundary conditions.
When the process reaches the $\varepsilon$-absorption layer it is projected
on the boundary $\partial D$ at a point $(x,y).$ The standard way
to replace the process inside the domain is based on a finite differences
approximations \cite{mascagni-simonov-2004}. The idea is to use a
normal approximation of the derivative at the boundary. In order to
simplify the description, we assume for instance that the boundary
is locally vertical that is it has the form $((x,z),y_{min}\leq z\leq y_{max})$.
Let $(x,y)$ be a point on this part of boundary, we have the order
one approximation 
\[
\frac{1}{2}\frac{u(x,y)-u(x+h,y)}{h}\simeq g(x,y)
\]
which leads to $u(x,y)\simeq2hg(x,y)+u(x+h,y).$ This simply means
in terms of randomization that $2hg(x,y)$ is added to the score of
the walk and that the process starts again in the domain at position
$(x+h,y).$ This approach does not give the elapsed time from position
$(x,y)$ to position $(x+h,y)$ and neither takes into account that
we deal with the Laplace operator. Another approach based on the randomization
of an integral equation has been introduced in \cite{simonov-1}.
This approach has an improved rate of convergence and can also deal
with transmission conditions. We propose now two alternative methods
one based on a kinetic approximation and the other on higher order
finite differences.

\subsubsection{Kinetic approximation}

A new scheme based on a kinetic approximation at the boundary has
been introduced in \cite{lejay-maire} for the simulation of multidimensional
diffusions in a media where the diffusion coefficient present some
discontinuities. The homogeneous Neumann boundary conditions were
also treated using this new scheme on a hard test case taken from
the couplex exercices. This scheme is based on a small parameter approximation
of the diffusion operator by a neutron transport operator. Another
application was the Monte Carlo solution of the Poisson-Boltzmann
equation in molecular dynamics \cite{bossy-etal-2010}. All the numerical
tests were very satisfactory compared to the finite differences approach.
The new scheme is proven to have an increased order of convergence
on the test case of a single sphere for the Poisson-Boltzmann equation.

We shall now describe this kinetic approximation in dimension 2 and
how to adapt it to the treatment of inhomogeneous Neumann boundary
conditions. We assume for the sake of simplicity that the process
has hitten from the right a vertical boundary with a Neumann condition
at a point $(x,y)$ which is infinitely away from the other boundaries.
We introduce a small parameter $h$ in order to approximate the Laplace
operator by the transport operator. We also pick a collision time
$t_{c}$ according to an exponential law of parameter 1 and a velocity
$(v_{x}=\cos\theta,v_{y}=\sin\theta)$ uniformly on the half unit
circle $(-\frac{\pi}{2}\leq\theta\leq\frac{\pi}{2}).$ The new position
is $(x+hv_{x}t_{c},y+hv_{y}t_{c})$ and the elapsed time is $h^{2}t_{c}.$
To deal with Neumann boundary conditions, we can write the Taylor
expansion 
\[
u(x+hv_{x}t_{c},y+hv_{y}t_{c})=u(x,y)+h\psi_{1}(x,y)+h^{2}\psi_{2}(x,y)+O(h^{3})
\]
 with 
\[
\psi_{1}(x,y)=-2v_{x}t_{c}g(x,y)+v_{y}t_{c}\frac{\partial u(x,y)}{\partial y}
\]
 and 
\[
\psi_{2}(x,y)=\frac{1}{2}\left(v_{x}^{2}t_{c}^{2}\frac{\partial^{2}u(x,y)}{\partial x^{2}}+v_{y}^{2}t_{c}^{2}\frac{\partial^{2}u(x,y)}{\partial y^{2}}+2v_{x}v_{y}t_{c}^{2}\frac{\partial^{2}u(x,y)}{\partial x\partial y}\right).
\]
Taking its mean value, we obtain 
\[
u(x,y)=\mathbb{E}(u(x+hv_{x}t_{c},y+hv_{y}t_{c}))+\frac{4g(x,y)}{\pi}h+f(x,y)h^{2}+O(h^{3})
\]
 because 
\[
\mathbb{E}\left(v_{y}\right)=0,\quad\quad\mathbb{E}\left(v_{x}\right)=\int_{-\frac{\pi}{2}}^{^{\frac{\pi}{2}}}\cos\theta d\theta,\quad\quad\mathbb{E}\left(t_{c}^{2}\right)=2.
\]
 Our scheme is the randomization of the previous formula: the motion
continues at point $(x+hv_{x}t_{c},y+hv_{y}t_{c})$ and the quantity
\[
\frac{4g(x,y)}{\pi}h+f(x,y)h^{2}
\]
 is added to the current score of the walk and $h^{2}$ to the total
time. The error is a $O(h^{3})$ each time the boundary is hitten
during the walk if we do not take into account the other boundaries.
If $(x,y)$ is close to another part of the boundary, the point 
\[
(x+hv_{x}t_{c},y+hv_{y}t_{c})
\]
 may be outside the domain. In this case, we reduce $h$ iteratively
by a factor 2 until it is inside. In the case of a general boundary,
the choice of the velocity law must be obviously adapted to the form
of the boundary to ensure that the process reenters in $D.$

\subsubsection{Order three finite differences approximation}

Some more sophisticated finite differences schemes can also be
considered. Given a step $h,$ we approximate the Laplace operator
and the normal derivative at the point $(x,y)$ of the boundary using
order 3 finite differences. For the Laplace operator, we have thanks
to the so-called order 2 diamond scheme 
\[
-\frac{1}{2}\Delta u(x,y)=f(x,y)\simeq\frac{-u(x-h,y)-u(x,y-h)-u(x,y+h)-u(x+h,y)+4u(x,y)}{2h^{2}}
\]
and for the centered normal derivative 
\[
\dfrac{1}{2}\frac{\partial u(x,y)}{\partial x}\simeq\frac{u(x+h,y)-u(x-h,y)}{4h}=g(x,y)+O(h^{2}).
\]
 By combining these two equations and getting rid of the fictitious
value $u(x-h,y),$ we obtain

\[
u(x,y)=\frac{2u(x+h,y)+u(x,y+h)+u(x,y-h)}{4}+\frac{h^{2}}{2}f(x,y)+2hg(x,y)+O(h^{3}).
\]
 This formula can be used for Monte Carlo simulations. Each time the
Brownian motion hits the boundary at a point $(x,y)$, the quantity
\[
\frac{h^{2}}{2}f(x,y)+2hg(x,y)
\]
 is added to the current score and the motion continues at one of
the positions 
\[
(x+h,y),(x,y+h),(x,y-h)
\]
 with the discrete probability law $(\frac{1}{2},\frac{1}{4},\frac{1}{4}).$
Some algorithms of a similar type for a walk on a grid have been introduced
in \cite{mikhailov-makarov-1997} and are called sliding on the boundary
methods because some of the possible points of replacement are on
the boundary. 

With our approach, if the boundary is not a straight line, the points
$(x,y+h)$ and $(x,y-h)$ may lie outside the domain. If it happens,
one needs to project the point on the boundary which induces an additional
error. We introduce now a new similar scheme for a general boundary
which avoids this projection. Instead of writing the previous equations
at point $(x,y),$ we write them at point $(x+h,y).$ We obtain 
\[
\frac{-u(x,y)-u(x+h,y-h)-u(x+h,y+h)-u(x+2h,y)+4u(x+h,y)}{2h^{2}}\simeq f(x+h,y)
\]
 for the Laplace operator. For the normal derivatives, we write the
two Taylor expansions 
\[
u(x+h,y)=u(x,y)+h\frac{\partial u(x,y)}{\partial x}+\frac{h^{2}}{2}\frac{\partial^{2}u(x,y)}{\partial x^{2}}+O(h^{3})
\]
 and 
\[
u(x+2h,y)=u(x,y)+2h\frac{\partial u(x,y)}{\partial x}+2h^{2}\frac{\partial^{2}u(x,y)}{\partial x^{2}}+O(h^{3})
\]
 which leads to 
\[
g(x,y)=\dfrac{1}{2}\frac{\partial u(x,y)}{\partial x}=\frac{-3u(x,y)-u(x+2h,y)+4u(x+h,y)}{4h}+O(h^{2}).
\]
 We finally obtain 
\[
u(x,y)=2g(x,y)h+f(x+h,y)h^{2}+\frac{u(x+h,y+h)+u(x+h,y-h)}{2}+O(h^{3})
\]
 which leads to the following Monte Carlo representation. The motion
is replaced equiprobably at point $(x+h,y+h)$ or at point $(x+h,y-h)$,
\[
2g(x,y)h+f(x+h,y)h^{2}
\]
 is added to the total score and $h^{2}$ to the total time. This
time, the motion is more likely to be in the domain. If not, we divide
$h$ iteratively by a factor 2 until it is, as we did for the kinetic
approximation.

\section{Numerical results \label{sec:Numerical-results}}

We shall now test our schemes on two Poisson equations in the square
$D=\left[-1,1\right]^{2}$ with increasing levels of difficulty. In
the first one, we consider a mixed Dirichlet-Neumann where there is
no problem of uniqueness for the solution and where the trajectories
stop when they hit the ``Dirichlet'' boundary. This enables a first
comparison between the different schemes and validate the new ones
introduced in section \ref{sec:WOS}. In the second one, we study
the pure Neumann problem for which additional difficulties arise,
especially the problems of the uniqueness of the solution and the
choice of the time $T_{0}$ to stop the walks.

\subsection{Mixed Dirichlet-Neumann Poisson equation}

Our first test case is the Poisson equation in the square domain $D=[-1,1]^{2}$
with boundary $\partial D=\partial D_{1}\cup\partial D_{2}$ defined
by
\[
\begin{cases}
\partial D_{1} & :=\left\{ (x,y),\left|y\right|=1,-1\leq x\leq1\right\} \cup\left\{ (x,y),x=-1,-1\leq y\leq1\right\} ,\\
\partial D_{2} & :=\left\{ (x,y),x=1,-1\leq y\leq1\right\} .
\end{cases}
\]
 We have a Neumann boundary condition on $\partial D_{1}$ and a Dirichlet
boundary condition on $\partial D_{2}$: 
\[
\begin{cases}
-\frac{1}{2}\Delta u(x,y) & =f(x,y):=-\alpha^{2}\exp(\alpha(x+y))\quad\text{{for}}\quad(x,y)\in D\\
\frac{1}{2}\frac{\partial u}{\partial n}(x,y) & =g_{1}(x,y):=\pm\frac{\alpha}{2}\exp(\alpha(x+y))\quad\text{{for}}\quad(x,y)\in\partial D_{1}\\
u(x,y) & =g_{2}(x,y):=\exp(\alpha(x+y))\quad\text{{for}}\quad(x,y)\in\partial D_{2},
\end{cases}
\]
where the sign of $g_{1}$ is negative for $y=+1$ and positive otherwise.
The stochastic process associated to this equation is a standard Brownian
motion $(B_{t},t\geq0)$, reflected on $\partial D_{1}$ and killed
on $\partial D_{2}.$ The solution $u$ has the stochastic representation
\[
u(x,y)=\mathbb{E}_{(x,y)}\left(\int_{0}^{\tau_{D}}f(B_{s})ds+\int_{0}^{\tau_{D}}g_{1}(B_{s})dV_{s}+g_{2}(B_{\tau_{D}})\right)=\exp(\alpha(x+y))
\]
 where $\tau_{D}$ is the first hitting time of $\partial D_{2}$.
The parameter $\alpha$ is introduced to obtain solutions with different
variations and consequently different variances in the Feynman-Kac
representations. We choose the values $\alpha_{1}=\frac{1}{3}$ ,
$\alpha_{2}=\frac{2}{3}$ and $\alpha_{3}=1$ ranked by increasing
degree of difficulty. We compute the solution using $N=50000$ Monte
Carlo simulations at some reference points $M_{1}=(0.8,0)$, $M_{2}=(0,0)$
and $M_{3}=(-0.8,0)$ which are located at positions at different
distances from $\partial D_{2}.$ We also compute the ratio $\frac{\sigma}{\sqrt{N}}$,
where $\sigma^{2}$ is the variance of the method which gives an estimation
of the Monte Carlo error.

\subsubsection{Euler Scheme approximations}

We denote by $e(\delta,\xi)$ the absolute error on the solution using
the reflected Euler scheme with the half space approximation on the
Dirichlet side. This error is computed in table \ref{tableau1DirichletNeumann}
with 3 different sets of parameters $(\delta_{1},\xi_{1})=(0.01,0.01)$
, $(\delta_{2},\xi_{2})=(0.001,0.000001)$ and $(\delta_{3},\xi_{3})=(0.001,0.001).$

\begin{table}
\centering{}%
\begin{tabular}{|c|c|c|c|c|c|c|c|}
\hline 
$(\alpha,M)$  & exact  & $e(\delta_{1},\xi_{1})$  & $\frac{\sigma_{1}}{\sqrt{N}}$  & $e(\delta_{2},\xi_{2})$  & $\frac{\sigma_{2}}{\sqrt{N}}$  & $e(\delta_{3},\xi_{3})$  & $\frac{\sigma_{3}}{\sqrt{N}}$\tabularnewline
\hline 
\hline 
$(\alpha_{1},M_{1})$  & 1.306  & 9.9E-4  & 2.0E-3  & 4.7E-3  & 1.7E-2  & 3.8E-4  & 2.0E-3\tabularnewline
\hline 
$(\alpha_{2},M_{1})$  & 1.705  & 1.8E-3  & 4.6E-3  & 4.6E-2  & 4.2E-2  & 2.1E-4  & 4.7E-3\tabularnewline
\hline 
$(\alpha_{3},M_{1})$  & 2.226  & 1.1E-2  & 8.9E-3  & 1.5E-1  & 8.9E-2  & 2.4E-3  & 9.0E-3\tabularnewline
\hline 
$(\alpha_{1},M_{2})$  & 1  & 1.6E-3  & 3.7E-3  & 6.1E-2  & 3.6E-2  & 3.2E-3  & 3.7E-3\tabularnewline
\hline 
$(\alpha_{2},M_{2})$  & 1  & 9.1E-3  & 8.2E-3  & 1.2E-2  & 8.4E-2  & 9.0E-3  & 8.3E-3\tabularnewline
\hline 
$(\alpha_{3},M_{2})$  & 1  & 2.6E-2  & 1.5E-2  & 6.4E-2  & 1.7E-1  & 1.8E-2  & 1.5E-2\tabularnewline
\hline 
$(\alpha_{1},M_{3})$  & 0.766  & 1.1E-2  & 4.1E-3  & 2.6E-2  & 4.3E-2  & 1.0E-3  & 4.2E-3\tabularnewline
\hline 
$(\alpha_{2},M_{3})$  & 0.587  & 1.6E-2  & 8.9E-3  & 4.7E-3  & 9.6E-2  & 4.7E-3  & 9.0E-3\tabularnewline
\hline 
$(\alpha_{3},M_{3})$  & 0.449  & 1.3E-2  & 1.6E-2  & 7.3E-2  & 1.9E-1  & 1.1E-2  & 1.7E-2\tabularnewline
\hline 
\end{tabular}\caption{Errors in the approximation of the solution of the mixed Dirichlet-Neumann
Poisson problem in the domain $D=[-1,1]^{2}$ with an Euler scheme
with parameters $(\delta_{1},\xi_{1})=(0.01,0.01)$ , $(\delta_{2},\xi_{2})=(0.001,0.000001)$
and $(\delta_{3},\xi_{3})=(0.001,0.001).$ \label{tableau1DirichletNeumann}}
\end{table}

For the parameters $(\delta_{1},\xi_{1}),$ we observe a good accuracy
of at least 2.6E-2 on the solution at all the reference points with
the various levels of difficulty. We note that the approximate solution
is more accurate and the variance smaller at point $M_{1}$ which
is very close to the Dirichlet boundary. The approximations at the
two other reference points are similar in terms of error and variance.
The CPU times on a standard PC for the computation of the solution
for the 3 levels of difficulty simultaneously are about 2 seconds
for $M_{1},$ 7 seconds for $M_{2}$ and 10 seconds for $M_{3}.$
Moreover, variance and bias increase with the level of difficulty.
To improve the accuracy of our method, we have first chosen to reduce
drastically the regularisation parameter and by a factor 10 the time discretization
parameter. The corresponding results for parameters $(\delta_{2},\xi_{2})$
were not satisfactory because the variance increases too much with
$\xi_{2}.$ This is especially true when $\alpha=\alpha_{3}$ where
the accuracy is no more than one digit. The results are a lot better
with $(\delta_{3},\xi_{3})$ with roughly the same variance than with
$(\delta_{1},\xi_{1})$ but with a smaller bias. For example, the
accuracy is 3 times better at point $M_{1}$ for the 3 values of $\alpha.$
The CPU times increase by a factor 10 which corresponds to the reduction
of $\delta.$

\subsubsection{Walk on spheres approximations\label{sub:WOS_approx}}

In table \ref{tableauWOSmixeddirichneum}, we compare three different
methods to compute the absolute errors on the exact solution all relying
on the walk on spheres method with absorption parameter $\varepsilon=10^{-6}$
but with different ways to handle the Neumann boundary conditions.
The first two errors $F_{3}(h)$ and $F_{2}(h)$ are based on the
finite differences method with scores respectively $g(x,y)h+f(x+h,y)h^{2}$
and $g(x,y)h$ at the boundary $\partial D_{1}.$ This enables to
emphasizes the differences between our new approach with the additional
term $f(x+h,y)h^{2}$ and the standard one. The last error $K(h)$
is computed thanks to the kinetic approximation. The simulations are
performed with two values $h_{1}=0.2$ and $h_{2}=0.1$ of the parameter
$h.$ The exit time of the unit circle and the associated uniform
position before absorption are pre-computed and stored in files of
size $10^{6}.$ 
\begin{table}
\begin{centering}
\begin{tabular}{|c|c|c|c|c|c|c|c|}
\hline 
$(\alpha,M)$  & $F_{3}(h_{1})$  & $F_{3}(h_{2})$  & $K(h_{1})$  & $K(h_{2})$  & $F_{2}(h_{1})$  & $F_{2}(h_{2})$  & $\frac{\sigma}{\sqrt{N}}$\tabularnewline
\hline 
$(\alpha_{1},M_{1})$  & 9.2E-4  & 2.4E-3  & 4.6E-3  & 1.6E-3  & 1.1E-2  & 4.4E-3  & 2E-3\tabularnewline
\hline 
$(\alpha_{2},M_{1})$  & 6.2E-3  & 5.3E-3  & 2.0E-2  & 1.3E-3  & 3.6E-2  & 1.9E-2  & 4E-3\tabularnewline
\hline 
$(\alpha_{3},M_{1})$  & 1.4E-2  & 7.7E-3  & 5.5E-2  & 3.5E-3  & 9.2E-2  & 4.9E-2  & 8E-3\tabularnewline
\hline 
$(\alpha_{1},M_{2})$  & 1.1E-2  & 2.3E-3  & 3.7E-3  & 2.8E-3  & 5.5E-2  & 2.4E-2  & 4E-3\tabularnewline
\hline 
$(\alpha_{2},M_{2})$  & 1.3E-2  & 1.8E-3  & 3.9E-2  & 3.3E-3  & 1.2E-1  & 8.5E-2  & 8E-3\tabularnewline
\hline 
$(\alpha_{3},M_{2})$  & 1.8E-2  & 9.0E-3  & 1.4E-1  & 2.5E-2  & 4.2E-1  & 2.0E-1  & 1E-2\tabularnewline
\hline 
$(\alpha_{1},M_{3})$  & 1.7E-2  & 1.1E-2  & 2.5E-2  & 1.2E-2  & 6.3E-2  & 4.2E-2  & 4E-3\tabularnewline
\hline 
$(\alpha_{2},M_{3})$  & 1.2E-2  & 1.0E-2  & 4.1E-2  & 1.1E-2  & 2.7E-1  & 1.4E-1  & 8E-3\tabularnewline
\hline 
$(\alpha_{3},M_{3})$  & 1.3E-2  & 6.2E-3  & 1.5E-2  & 3.8E-3  & 5.8E-1  & 3.1E-1  & 2E-2\tabularnewline
\hline 
\end{tabular}
\par\end{centering}

\caption{Errors in the approximation of the solution of the mixed Dirichlet-Neumann
Poisson problem in the domain $D=[-1,1]^{2}$ with WOS approximations.
The errors $F_{3}(h)$ and $F_{2}(h)$ are based on the finite differences
method with scores respectively $g(x,y)h+f(x+h,y)h^{2}$ and $g(x,y)h$
at the boundary $\partial D_{2}$. The error $K(h)$ is computed thanks
to the kinetic approximation. $h_{1}=0.2$ and $h_{2}=0.1$.\label{tableauWOSmixeddirichneum}}
\end{table}

The CPU times for the computation of the solution at points $(M_{1},M_{2},M_{3})$
with discretization step $h_{1}$ are $(4,14,21)$ for the finite
differences method, and $(5,20,30)$ for the kinetic approximation
(they are twice bigger when using $h_{2}).$ In terms of accuracy,
we observe that the errors $F_{2}(h)$ are clearly a lot bigger than
the ones obtained with the two new methods. The errors $F_{3}(h_{2})$
and $K(h_{2})$ are similar and are furthermore comparable to $e(\delta_{3},\xi_{3})$
for a CPU time twice smaller. The error $F_{3}(h_{1})$ is similar
to $e(\delta_{1},\xi_{1})$ for a CPU time twice bigger. We can conclude
that the two new schemes are very efficient especially when one desires
an accurate solution. Nevertheless and unlike with Dirichlet boundary
conditions, there is not such a big difference in terms of efficiency
between Euler schemes and WOS methods.

\subsection{Pure Neumann Poisson Equation}

\subsubsection{Preliminary example}

\begin{figure}
\centering{}\includegraphics[scale=0.85]{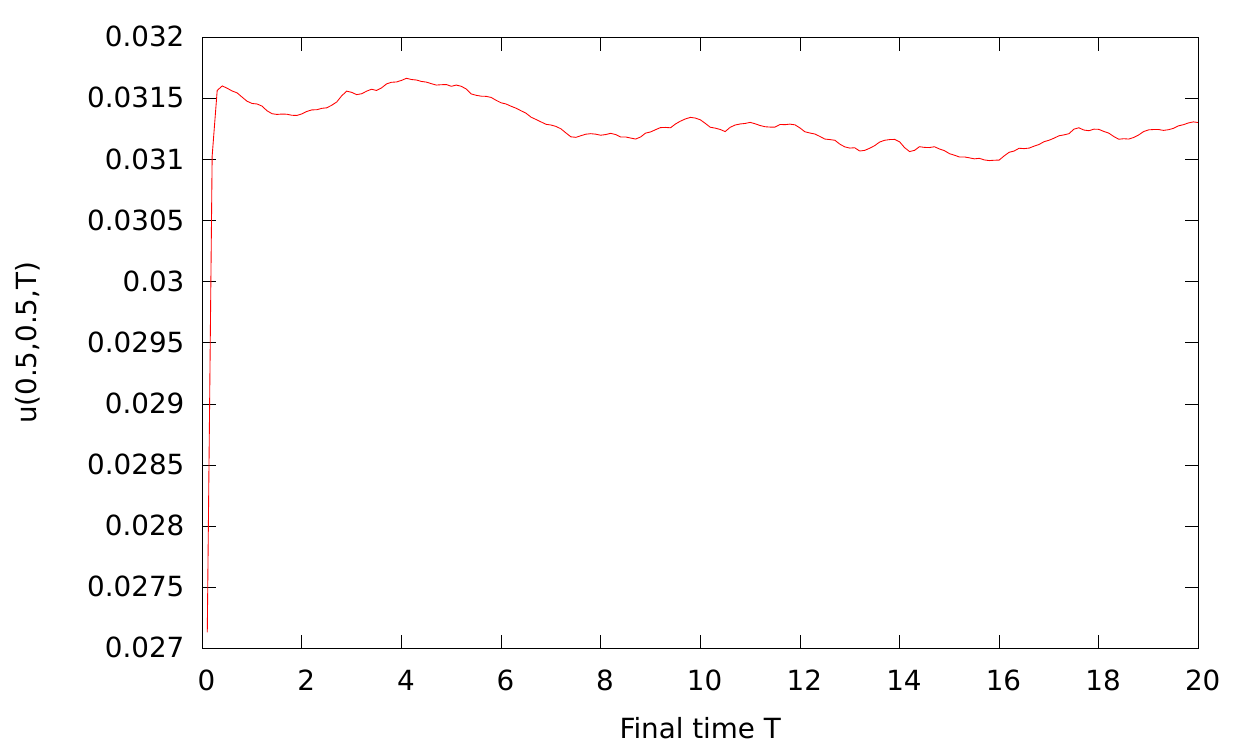}\caption{Estimation of the solution of the pure Neumann problem at point $(x,y)=(-0.5,-0.5)$
in function of the final time $T$ \label{fig:solution0505}}
\end{figure}

In this part, we illustrate our theoretical results of section \ref{sec:properties_stochastic_representation}
on the solution of a very simple Cauchy problem and of its related
variance as a function of $T.$ We consider the Poisson equation

\[
\frac{1}{2}\Delta u=-f(x,y)=-2(3x^{2}-1)(y^{2}-1)^{2}-2(3y^{2}-1)(x^{2}-1)^{2}
\]
 for $(x,y)\in D=[-1,1]^{2}$ and homogeneous Neumann boundary conditions
$\frac{\partial u}{\partial n}=0$ on $\partial D$. The stochastic
process associated to this equation is a standard reflected Brownian
motion and its invariant measure is the uniform law in $D.$ The exact
solution with mean value zero is 
\[
u(x,y)=(x^{2}-1)^{2}(y^{2}-1)^{2}-\frac{64}{225}.
\]
 In this simple domain, we know that the second leading eigenvalue
$\lambda_{1}$ of the operator $-\frac{1}{2}\Delta$ with
pure Neumann boundary condition is $\frac{\pi^{2}}{8}$. This means
that the convergence of the solution of the Cauchy problem towards
the solution of the above equation is a $O(\exp(-\frac{\pi^{2}}{8}T)).$
We have also proven that the main part of the variance of our scheme
increases linearly as $C_{3}T$, where 
\[C_{3}=\frac{1}{4}\int{}_{D}\left|\nabla u\right|^{2}.
\]
Here, we have 
\[
\left|\nabla u\right|^{2}=(4x(x^{2}-1)(y^{2}-1)^{2})^{2}+(4y(y^{2}-1)(x^{2}-1)^{2})^{2}
\]
 and finally, 
 \[
 C_{3}=\frac{32768}{33075}\simeq0.99.
 \] 
 We compute an
approximate solution $u(x,y,T)$ at point $(x,y)=(-0.5,-0.5)$ and
its related variance for values of $T\in[0.1,20].$ The numerical
method used is the Euler scheme with a small stepsize $\delta=0.001$
and a huge number $N=5\times10^{7}$ of simulations. The exact solution
with mean value zero is 
\[
u(-0.5,-0.5)=\frac{81}{256}-\frac{64}{225}\simeq0,03196.
\]
 In figure \ref{fig:solution0505}, we observe that $u(0.5,0.5,T)$
converges quickly to a constant (modulo some statistical variations)
which is close to $u(-0.5,-0.5).$

\begin{figure}
\centering{}\includegraphics{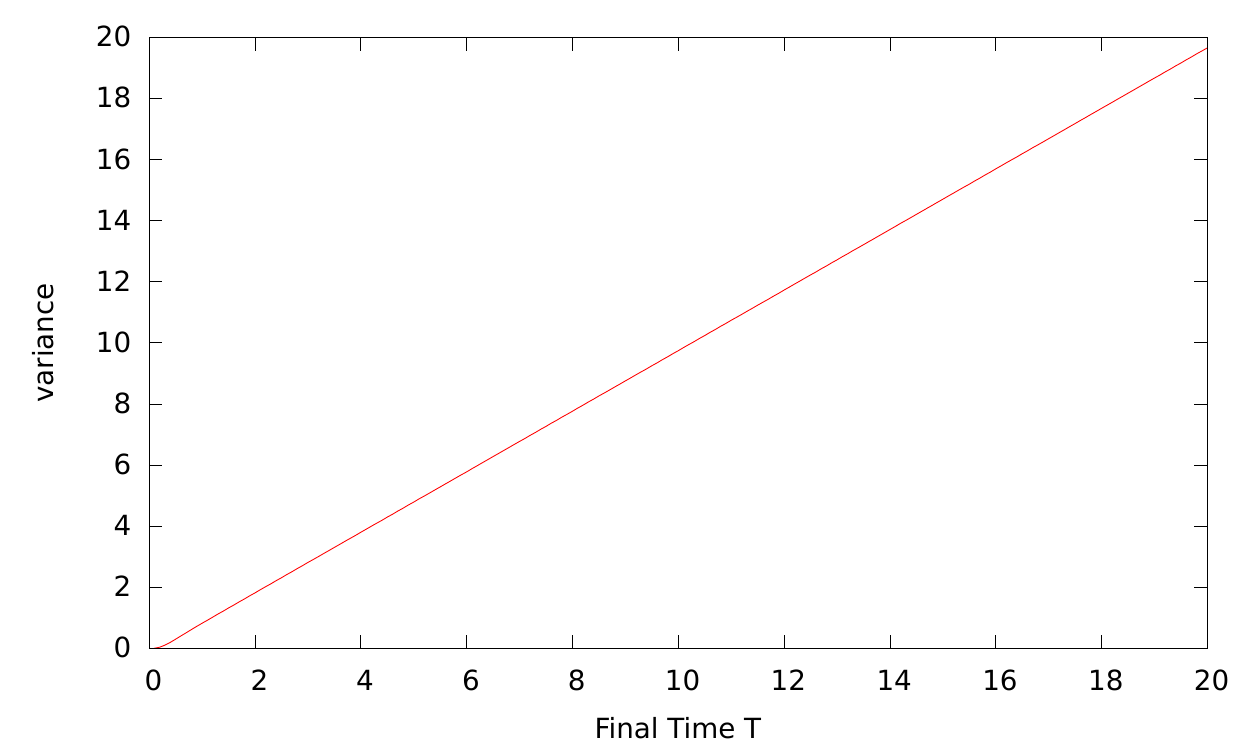}\caption{Estimation of the variance of the approximation of the solution of
the pure Neumann problem at point $(x,y)=(-0.5,-0.5)$ in function
of the final time $T$\label{fig:variance0505}}
\end{figure}

Figure \ref{fig:variance0505} concerns the variance: we observe that
it increases linearly as a function of $T.$ If we compute the slope
of the variance using for instance a linear regression on the approximate
values at times $T_{i}=8+i$ for $1\leq i\leq8,$ we obtain about
$0.99.$ This means that the non linear part which behaves as $C_{2}\sqrt{T}$
should be very small on this particular example.

On this simple example, we have been able to confirm the results obtained
in section \ref{sec:properties_stochastic_representation} in the
more general setting of a domain in dimension 2. However, we should
now consider inhomogeneous Neumann boundary conditions where both
bias and variance may increase due to local time approximation.

\subsubsection{Numerical parameters and error criteria}

Our main test case is the Poisson equation in the square domain $D=[-1,1]^{2}$
defined by 
\[
-\frac{1}{2}\Delta u(x,y)=-\alpha^{2}\exp(\alpha(x+y))
\]
 in $D$ and with Neumann boundary conditions $\frac{1}{2}\frac{\partial u}{\partial n}=\pm\frac{\alpha}{2}\exp(\alpha(x+y))$
on $\partial D$, where the sign is positive on the bottom and left
sides of the boundary and negative on the right and top sides. The
solution of this equation with zero mean value with respect to the
invariant measure is hence 
\[
u(x,y)=\exp(\alpha(x+y))-\frac{1}{4}\int_{-1}^{1}\int_{-1}^{1}\exp(\alpha(x+y))dxdy
\]
 that is 
\[
u(x,y)=\exp(\alpha(x+y))-\frac{(\exp(\alpha)-\exp(-\alpha))^{2}}{4\alpha^{2}}.
\]
 This solution is the one that we are likely to obtain numerically
if we have a perfect simulation of the reflected Brownian motion and
a good choice of the time $T_{0}$ when we stop the trajectories.
Even though we have noticed on our test cases on the mixed problem
that our numerical schemes are quite efficient, they are obviously
not perfect and consequently introduce supplementary errors. In fact,
regardless to the usual Monte Carlo and discretization errors, we do
not compute the solution of the equation with zero mean value with
respect to the invariant measure but another one. Nevertheless, we
know that our approximation $\bar{u}(x,y)$ should be of the form
$u(x,y)+a$ where $a$ is a constant. In order to check the quality
of our approximation and estimate this constant, we will compute the
minimum of the weighted cost function 
\[
J(a)=\frac{1}{\pi^{2}}\int_{-1}^{1}\int_{-1}^{1}\frac{(\bar{u}(x,y)-u(x,y)-a)^{2}}{\sqrt{(1-x^{2})(1-y^{2})}}dxdy.
\]
 Using a Gauss quadrature formula on the Tchebychef grid, $J(a)$
can be approximated by 
\[
J_{1}(a)=\frac{1}{P^{2}}\sum_{i=1}^{P}\sum_{j=1}^{P}(\bar{u}(x_{i},x_{j})-u(x_{i},x_{j})-a)^{2}
\]
 using the Tchebychef points $x_{i}=\cos(\frac{(2i-1)\pi}{2P}).$
The minimum of $J_{1}$ is achieved for 
\[
\overline{a}=\frac{1}{P^{2}}\sum_{i=1}^{P}\sum_{j=1}^{P}(\bar{u}(x_{i},x_{j})-u(x_{i},x_{j}))
\]
 which indicates the bias with respect to the perfectly simulated
solution while the value $\rho=\sqrt{J_{1}(\overline{a})}$ quantifies
the adequacy to the model. The choice of the weighted cost function
is motivated by the use of bidimensional Tchebychef interpolation
polynomials in the stochastic spectral methods of the next section.
Indeed, they rely on the pointwise approximations at the points of
the Tchebychef grid. In practice, we choose $P=3$ which is sufficient
to have a good accuracy on the weighted integrals. Once $\overline{a}$
and $\rho$ have been computed, we compute the approximate solution
$F(M)=\bar{u}(x,y)$ and the error on the model $B(M)=\left|\bar{u}(x,y)-u(x,y)-\overline{a}\right|$
for two sets of points $M.$ The first one contains points $M_{4}=(0,0)$,
$M_{5}=(-0.2,0.2)$, $M_{6}=(-0.8,0.8)$ and the second one contains
points $M_{7}=(0,0.8)$, $M_{8}=(0.2,0.6)$, $M_{9}=(0.4,0.4).$ Morevover
all computations are performed using $N=50000$ Monte Carlo simulations.

\subsubsection{Euler scheme approximations}

First of all, we shall choose $T_{0}=10,$ which corresponds to a bias equal
to $\exp(-1.25\pi^{2})\simeq5\times10^{-6}$ in all the following
numerical tests. This value is small enough so that this bias is negligible
with respect to the other errors and the variance not too large.

The values of the quantities $\frac{\sigma}{\sqrt{N}}$ depend weakly
on the starting point and on the discretization. They are approximatively
equal to respectively 0.006, 0.013 and 0.03 for respectively $\alpha_{1}=\frac{1}{3},\,\alpha_{2}=\frac{2}{3}$
and $\alpha_{3}=1$ which is about 1.5 times bigger than in the mixed
Dirichlet-Neumann case. When $(\delta_{1},\xi_{1})=(0.01,0.01),$
$\overline{a}$ and $\rho$ are equal respectively to $-0.0059$,
$-0.052$, $-0.130$ and $0.0039$, $0.027$, $0.0498$. When $(\delta_{2},\xi_{2})=(0.001,0.001),$
$\overline{a}$ and $\rho$ are equal respectively to $-0.0015$,$-0.0042$,
$-0.016$ and $0.0053$,$0.021$, $0.042$.

We observe that the values of $\overline{a}$ and $\rho$ are small
which indicates a good adequacy to the model. Moreover, $\rho$ and
especially $\overline{a}$ are significantly closer to zero when the
discretisation parameters decrease from $(\delta_{1},\xi_{1})$ to
$(\delta_{2},\xi_{2}).$ This shows that the approximate solution
gets closer to $u(x,y)$ as it has been proven in section \ref{sec:properties_stochastic_representation}.
\begin{table}[th]
\centering{}%
\begin{tabular}{|c|c|c|c|c|c|c|}
\hline 
 & $F(M_{4})$  & $F(M_{5})$  & $F(M_{6})$  & $B(M_{4})$  & $B(M_{5})$  & $B(M_{6})$\tabularnewline
\hline 
\hline 
$(\alpha_{1},(\delta_{1},\xi_{1}))$  & $-0.045$  & $-0.042$  & $-0.042$  & $0.002$  & $0.001$  & $0.001$\tabularnewline
\hline 
$(\alpha_{2},(\delta_{1},\xi_{1}))$  & $-0.204$  & $-0.172$  & $-0.210$  & $0.005$  & $0.038$  & $0.001$\tabularnewline
\hline 
$(\alpha_{3},(\delta_{1},\xi_{1}))$  & $-0.480$  & $-0.506$  & $-0.508$  & $0.031$  & $0.005$  & $0.003$\tabularnewline
\hline 
$(\alpha_{1},(\delta_{2},\xi_{2}))$  & $-0.030$  & $-0.041$  & $-0.036$  & $0.01$  & $0.002$  & $0.003$\tabularnewline
\hline 
$(\alpha_{2},(\delta_{2},\xi_{2}))$  & $-0.155$  & $-0.182$  & $-0.146$  & $0.006$  & $0.02$  & $0.016$\tabularnewline
\hline 
$(\alpha_{3},(\delta_{2},\xi_{2}))$  & $-0.374$  & $-0.365$  & $-0.395$  & $0.023$  & $0.032$  & $0.002$\tabularnewline
\hline 
\end{tabular}\caption{Approximation $F$ of the unbiased solution $u_{\overline{a}}$ of
the pure Neumann problem with an Euler scheme with steps $(\delta_{1},\xi_{1})=(0.01,0.01)$
and $(\delta_{2},\xi_{2})=(0.001,0.001)$ at points $M_{4}=(0,0)$,
$M_{5}=(-0.2,0.2)$, $M_{6}=(-0.8,0.8)$ and the corresponding errors
$B(M)=\left|\bar{u}(M)-u(M)-\overline{a}\right|$.\label{tab:approximationnonbiaisee1} }
\end{table}

On this first set of points, we can see that for the same value of
$(\alpha,\delta,\xi)$ the direct estimations $F(M_{4}),\, F(M_{5})$
and $F(M_{6})$ are close to each others. We observe also that the
approximation model plays an important role. Indeed, the direct approximations
are significantly different from each other for the two sets of parameters
but nevertheless the errors $B(M)$ are quite small for both sets.
The maximum absolute errors are 0.01, 0.038 and 0.031 for respectively
$\alpha_{1},\,\alpha_{2}$ and $\alpha_{3}$. 
\begin{table}
\centering{}%
\begin{tabular}{|c|c|c|c|c|c|c|}
\hline 
 & $F(M_{7})$  & $F(M_{8})$  & $F(M_{9})$  & $B(M_{7})$  & $B(M_{8})$  & $B(M_{9})$\tabularnewline
\hline 
$(\alpha_{1},(\delta_{1},\xi_{1}))$  & $0.266$  & $0.255$  & $0.249$  & $0.004$  & $0.007$  & $0.013$\tabularnewline
\hline 
$(\alpha_{2},(\delta_{1},\xi_{1}))$  & $0.550$  & $0.524$  & $0.526$  & $0.054$  & $0.028$  & $0.031$\tabularnewline
\hline 
$(\alpha_{3},(\delta_{1},\xi_{1}))$  & $0.681$  & $0.715$  & $0.692$  & $0.033$  & $0.001$  & $0.022$\tabularnewline
\hline 
$(\alpha_{1},(\delta_{2},\xi_{2}))$  & $0.249$  & $0.272$  & $0.271$  & $0.02$  & $0.005$  & $0.004$\tabularnewline
\hline 
$(\alpha_{2},(\delta_{2},\xi_{2}))$  & $0.549$  & $0.537$  & $0.518$  & $0.006$  & $0.006$  & $0.025$\tabularnewline
\hline 
$(\alpha_{3},(\delta_{2},\xi_{2}))$  & $0.832$  & $0.897$  & $0.790$  & $0.004$  & $0.068$  & $0.038$\tabularnewline
\hline 
\end{tabular}\caption{Approximation $F$ of the unbiased solution $u_{\overline{a}}$ of
the pure Neumann problem with an Euler scheme with steps $(\delta_{1},\xi_{1})=(0.01,0.01)$
and $(\delta_{2},\xi_{2})=(0.001,0.001)$ at points $M_{7}=(0,0.8)$,
$M_{8}=(0.2,0.6)$, $M_{9}=(0.4,0.4)$ and the corresponding errors
$B(M)=\left|\bar{u}(M)-u(M)-\overline{a}\right|$.}
\end{table}

The same conclusions hold for the second set of points. The maximum
absolute errors are slightly bigger 0.02, 0.054 and 0.068 for respectively
$\alpha_{1},\,\alpha_{2}$ and $\alpha_{3}.$ These maximum errors
are at least 2 or 3 times bigger than in the mixed Dirichlet-Neumann
case. Furthermore the CPU times are about twice larger in the pure
Neumann Case. We can conclude that the pure Neumann problem is a lot
harder to solve but that our algorithm still provides an acceptable
accuracy for pointwise approximations.

\subsubsection{Walk on spheres approximations}

We have noticed in section \ref{sub:WOS_approx} that the two new
methods to handle the boundary conditions have the same accuracy.
We have chosen to use the order three finite differences in the following.
Furthermore to make our simulations, we have precomputed and stored
100 positions of $100000$ discretized trajectories as described in
section \ref{sub:WOS-pure-neumann}. The time to open this file is
negligible compared to the rest of the simulation times. They are
about twice larger than in the mixed Dirichlet-Neumann case. The quantities
$\frac{\sigma}{\sqrt{N}}$ are still approximatively equal to respectively
$0.006$, $0.013$ and $0.03$ for respectively $\alpha_{1}=\frac{1}{3},\,\alpha_{2}=\frac{2}{3}$
and $\alpha_{3}=1.$ For $h_{1}=0.1,$ $\overline{a}$ and $\rho$
are equal respectively to $-0.021$, $-0.106$, $-0.302$ and $0.018$,
$0.053$, $0.134$. For $h_{2}=0.05,$ $\overline{a}$ and $\rho$
are equal respectively to $-0.012$, $-0.060$, $-0.176$ and $0.012$,
$0.034$, $0.085$. We observe that the value of $\rho$ and especially
$\overline{a}$ are larger than with the Euler Scheme method.

\begin{table}
\begin{centering}
\begin{tabular}{|c|c|c|c|c|c|c|}
\hline 
 & $F(M_{4})$  & $F(M_{5})$  & $F(M_{6})$  & $B(M_{4})$  & $B(M_{5})$  & $B(M_{6})$\tabularnewline
\hline 
$(\alpha_{1},h_{1})$  & $-0.055$  & $-0.059$  & $-0.062$  & $0.002$  & $0.001$  & $0.004$\tabularnewline
\hline 
$(\alpha_{2},h_{1})$  & $-0.247$  & $-0.259$  & $-0.263$  & $0.016$  & $0.004$  & $0.001$\tabularnewline
\hline 
$(\alpha_{3},h_{1})$  & $-0.634$  & $-0.669$  & $-0.668$  & $0.049$  & $0.014$  & $0.015$\tabularnewline
\hline 
$(\alpha_{1},h_{2})$  & $-0.043$  & $-0.050$  & $-0.048$  & $0.006$  & $0.001$  & $0.001$\tabularnewline
\hline 
$(\alpha_{2},h_{2})$  & $-0.198$  & $-0.215$  & $-0.217$  & $0.019$  & $0.002$  & $0.001$\tabularnewline
\hline 
$(\alpha_{3},h_{2})$  & $-0.517$  & $-0.547$  & $-0.563$  & $0.04$  & $0.01$  & $0.005$\tabularnewline
\hline 
\end{tabular}
\par\end{centering}

\caption{Approximation $F$ of the unbiased solution $u_{\overline{a}}$ of
the pure Neumann problem with the WOS method at points $M_{4}=(0,0)$,
$M_{5}=(-0.2,0.2)$, $M_{6}=(-0.8,0.8)$ and the corresponding errors
$B(M)=\left|\bar{u}(M)-u(M)-\overline{a}\right|$. The parameters
are $h_{1}=0.1$ and $h_{2}=0.05$ }
\end{table}

Once again, the direct estimations $F(M_{4}),\, F(M_{5})$ and $F(M_{6})$
are close to each others for a given value of the parameters $(\alpha,h)$.
The maximum absolute errors are 0.006, 0.016 and 0.049 for respectively
$\alpha_{1},\,\alpha_{2}$ and $\alpha_{3}.$ This shows that even
if the solution computed is further away from $u(x,y),$ the accuracy
is similar than the one obtained with the Euler scheme method. 
\begin{table}[th]
\centering{ %
\begin{tabular}{|c|c|c|c|c|c|c|}
\hline 
 & $F(M_{7})$  & $F(M_{8})$  & $F(M_{9})$  & $B(M_{7})$  & $B(M_{8})$  & $B(M_{9})$\tabularnewline
\hline 
$(\alpha_{1},h_{1})$  & $0.231$  & $0.236$  & $0.240$  & $0.016$  & $0.011$  & $0.007$\tabularnewline
\hline 
$(\alpha_{2},h_{1})$  & $0.399$  & $0.412$  & $0.424$  & $0.043$  & $0.03$  & $0.018$\tabularnewline
\hline 
$(\alpha_{3},h_{1})$  & $0.455$  & $0.479$  & $0.503$  & $0.087$  & $0.063$  & $0.039$\tabularnewline
\hline 
$(\alpha_{1},h_{2})$  & $0.246$  & $0.253$  & $0.253$  & $0.011$  & $0.003$  & $0.003$\tabularnewline
\hline 
$(\alpha_{2},h_{2})$  & $0.466$  & $0.475$  & $0.474$  & $0.021$  & $0.012$  & $0.013$\tabularnewline
\hline 
$(\alpha_{3},h_{2})$  & $0.627$  & $0.623$  & $0.626$  & $0.041$  & $0.045$  & $0.042$\tabularnewline
\hline 
\end{tabular}} \caption{Approximation $F$ of the unbiased solution $u_{\overline{a}}$ of
the pure Neumann problem with the WOS method at points $M_{7}=(0,0.8)$,
$M_{8}=(0.2,0.6)$, $M_{9}=(0.4,0.4)$ and the corresponding errors
$B(M)=\left|\overline{u}(M)-u(M)-\overline{a}\right|$. The parameters
are $h_{1}=0.1$ and $h_{2}=0.05$ }
\end{table}

The maximum absolute errors are bigger for this set of points especially
for $(\alpha_{3},h_{1}).$ Nevertheless, the maximum absolute errors
are 0.011, 0.021 and 0.045 for respectively $\alpha_{1},\,\alpha_{2}$
and $\alpha_{3}.$ We can conclude that we achieve a good accuracy
on this pure Neumann problem but with an increased computational cost
compared to the mixed Dirichlet-Neumann problem.

\section{Stochastic spectral methods \label{sec:Stochastic-spectral-methods}}

\subsection{Spectral formulation}

In this section we describe how to adapt the stochastic spectral formulations
introduced in \cite{maire-tanre-2008} and studied in detail in \cite{maire-tanre-2009}
to the case of pure Neumann boundary conditions. These formulations
are similar to usual spectral methods based on polynomial approximations
\cite{canuto_etal-1988} but they are built using relevant information,
not necessarily at the collocation points, given by the Feynman-Kac
formula. They are an extension of the sequential Monte Carlo algorithms
for solving linear partial differential equations developed in \cite{gobet-maire-2004,gobet-maire-2005}.
These stochastic spectral formulations are asymptotically perfectly
conditioned and quite easy to build for Dirichlet boundary conditions..
The case of mixed boundary conditions is nor described nor studied
here because it is a straightforward extension of our previous works.
For pure Neumann conditions the situation is quite different because
of the non-uniqueness of the solution.We have to build basis functions
verifying centering conditions in order to obtain an invertible spectral
formulation. These new centering procedures either exact or approximate
are described in section 6.2.

When solving the pure Neumann problem using usual deterministic methods
like finite elements, one also has to take into account very accurately
the compatibility conditions and the non-uniqueness of the solution.
Two approaches are usually used. The first one consists in fixing
the value of the solution at a specified node in order to avoid the
resolution of a singular linear system. The second one leads to a
singular system but it is solved using an iterative method like the
conjuguate gradient for positive semi-definite linear system. In this
second case, it is extremely important that the compatibility condition
is verified at the discrete level which is obtained via discrete projectors.
All these questions as well as an accurate study of the condition
number are treated in \cite{bochev-lehoucq}.

Our stochastic formulation consists in computing a global linear approximation
$P_{N}u$ of the solution $u$ using its values at some points $x_{i}\in D.$
This global approximation writes 
\[
P_{N}u(x)=\sum_{i=1}^{N}u(x_{i})\Psi_{i}(x)
\]
 for some functions $\Psi_{i}(x)$ that are at least twice continuously
differentiable and we assume that they verify the centering condition
$\int_{D}\Psi_{i}(x)p(x)dx=0.$ Note that this last condition implies
that $P_{N}u$ belongs to the space of functions that have a zero
mean value with respect to $p(x)$ and ensures the uniqueness of the
solution. We also assume that for every point $x_{i}$, we can approximate
$u(x_{i})$ via for instance a numerical approximation of the Feynman-Kac
formula by 
\[
\widetilde{u}(x_{i})=\sum_{j=1}^{q}\alpha_{i,j}g(z_{i,j}^{b,\delta,T})+\sum_{j=1}^{p}\beta_{i,j}f(z_{i,j}^{s,\delta,T}).
\]
 In practice this approximation is also such that 
\[
\lim_{p,q,T\rightarrow\infty,\delta\rightarrow0}\widetilde{u}(x_{i})=u(x_{i})-\int_{D}u(x)p(x)dx
\]
 which indicates that the solution that we compute numerically is
close to the one with zero mean value with respect to $p(x)$. 
The
coefficients $\alpha_{i,j}$ and $\beta_{i,j}$ are positive weights.
The $q$ points $z_{i,1}^{b,\delta,T}, \cdots,z_{i,q}^{b,\delta,T}$ are located on the boundary
$\partial D$, the $p$ points $z_{i,1}^{s,\delta,T}, \cdots, z_{i,p}^{s,\delta,T}$ in $D$, $\delta$
stands for the discretization parameter of the simulated reflected
diffusion and $T$ is the deterministic time when we decide to stop
our random walk. We now let $r_{N}(x)=u(x)-P_{N}u(x)$ and write the
partial differential equation solved by $r_{N}(x).$ We have 
\[
Lr_{N}=Lu-LP_{N}u=- f - LP_{N}u
\]
 in $D$ with boundary conditions 
\[
\frac{\partial r_{N}}{\partial n_{a}}=g-\frac{\partial P_{N}u}{\partial n_{a}}
\]
 and hence the approximation 
\[
r_{N}(x_{i})=\sum_{j=1}^{q}\alpha_{i,j}\left(g(z_{i,j}^{b,\delta,T})-\frac{\partial P_{N}u(z_{i,j}^{b,\delta,T})}{\partial n_{a}}\right)+\sum_{j=1}^{p}\beta_{i,j}\left(f(z_{i,j}^{s,\delta,T}) +LP_{N}u(z_{i,j}^{s,\delta,T})\right)
\]
 which leads to the linear system $Cu=d$ with 
\[
c_{i,i}=\sum_{j=1}^{q}\alpha_{i,j}\frac{\partial\Psi_{i}(z_{i,j}^{b,\delta,T})}{\partial n_{a}} - \sum_{j=1}^{p}\beta_{i,j}L\Psi_{i}(z_{i,j}^{s,\delta,T})+1-\Psi_{i}(x_{i}),
\]
 for $i\neq k,$ 
\[
c_{i,k}=\sum_{j=1}^{q}\alpha_{i,j}\frac{\partial\Psi_{k}(z_{i,j}^{b,\delta,T})}{\partial n_{a}} - \sum_{j=1}^{p}\beta_{i,j}L\Psi_{k}(z_{i,j}^{s,\delta,T})-\Psi_{k}(x_{i})
\]
 and 
\[
d_{i}=\sum_{j=1}^{q}\alpha_{i,j}g(z_{i,j}^{b,\delta,T})+\sum_{j=1}^{p}\beta_{i,j}f(z^{s,\delta,T}_{i,j}).
\]
 As we have done in \cite{maire-tanre-2008,maire-tanre-2009}, we
can look at the asymptotic system we obtain when $p,q,T\rightarrow\infty$
and when $\delta\rightarrow0.$ The term 
\[
\sum_{j=1}^{q}\alpha_{i,j}\frac{\partial\Psi_{k}(z_{i,j}^{b,\delta,T})}{\partial n_{a}}-\sum_{j=1}^{p}\beta_{i,j}L\Psi_{k}(z_{i,j}^{s,\delta,T})
\]
 is our Monte Carlo approximation at point $x_{i}$ of the solution
of the equation 
\[
Lv=L\Psi_{k}
\]
 with boundary conditions 
\[
\frac{\partial v}{\partial n_{a}}=\frac{\partial\Psi_{k}}{\partial n_{a}}
\]
 on $\partial D$ that is $\Psi_{k}(x_{i})$ because $\int_{D}\Psi_{k}(x)p(x)dx=0.$
As in our previous papers, this shows immediately that the matrix
of the asymptotic system converges toward the identity matrix of size
$N.$ This also means that the condition number of the system is naturally
close to one even without additional preconditioning techniques like
for instance the ones developed in \cite{tang-vuik}.

\subsection{Centering procedures}

\subsubsection{Exact centering}

We now describe how this method works in practice. The main problem
is that in general, usual linear approximations do not verify the
centering conditions. We consider for the moment that we start with
$N+1$ Lagrange interpolation polynomials $\varphi_{i}$ at $N+1$
points $x_{i}\in D$. Such functions verify $\varphi_{i}(x_{k})=\delta_{i,k}.$
The usual polynomial interpolation $Q_{N+1}u$ of degree $N$ of a
function $u$ writes $Q_{N+1}u(x)=\sum_{i=1}^{N+1}u(x_{i})\varphi_{i}(x)$.
Considering the constant function $u=1,$we obviously have $\sum_{i=1}^{N+1}\varphi_{i}(x)=1$
and hence 
\[
\sum_{i=1}^{N+1}\int_{D}\varphi_{i}(x)p(x)dx=1,
\]
 which proves that this usual interpolation cannot verify the centering
conditions. Nevertheless, we can choose an index $i_{0}$ such that
\[
\left|\int_{D}\varphi_{i_{0}}(x)p(x)dx\right|=\max_{1\leq i\leq N+1}\left|\int_{D}\varphi_{i}(x)p(x)dx\right|\neq0
\]
 for a sake of stability in the following approximation. The centering
condition for $Q_{N+1}u$ writes 
\[
\sum_{i=1}^{N+1}\left(\int_{D}\varphi_{i}(x)p(x)dx\right)u(x_{i})=0
\]
 which leads to 
\[
Q_{N+1}u(x)=\sum_{i\neq i_{0}}u(x_{i})\left(\varphi_{i}(x)-\frac{\int_{D}\varphi_{i}(x)p(x)dx}{\int_{D}\varphi_{i_{0}}(x)p(x)dx}\varphi_{i_{0}}(x)\right).
\]
 Letting now 
 \[
 \Psi_{i}(x)=\varphi_{i}(x)-\frac{\int_{D}\varphi_{i}(x)p(x)dx}{\int_{D}\varphi_{i_{0}}(x)p(x)dx}\varphi_{i_{0}}(x),
 \]
the new basis functions verify the centering conditions $\int_{D}\Psi_{i}(x)p(x)dx=0$
and still $\Psi_{i}(x_{k})=\delta_{i,k}.$ This centering procedure
can be easily extended to general linear approximations not necessarily
of interpolation type. It is essentially the same approach than the
projection method described in \cite{bochev-lehoucq} for the finite
elements method.

\subsubsection{Approximate centering}

This solves our problem whenever the integrals $\int_{D}\varphi_{i}(x)p(x)dx$
can be computed exactly. This may happen when the domain is simple
and when the density $p(x)$ is known. This include for example the
case of the Poisson equation in a hypercube using an interpolation
on a Tchebychef grid. In many other cases, these integrals need to
be computed numerically. The density $p(x)$ can be approximated by
the law of the position $Y$ of $M$ particles moving according to
the reflected diffusion starting at any given point in $D$ at a time
$T_{1}$ large enough. We obtain 
\[
\int_{D}\varphi_{i}(x)p(x)dx\simeq\frac{1}{M}\sum_{j=1}^{M}\varphi_{i}(Y_{j})
\]
 and new basis functions which are defined by 
\[
\Psi_{i}(x)=\varphi_{i}(x)-\frac{\frac{1}{M}\sum_{j=1}^{M}\varphi_{i}(Y_{j})}{\frac{1}{M}\sum_{j=1}^{M}\varphi_{i_{0}}(Y_{j})}\varphi_{i_{0}}(x).
\]
 The coefficients of the asymptotic spectral matrix are
\[
\begin{cases}
c_{i,i} & =1-\int_{D}\Psi_{i}(x)p(x)dx\\
c_{i,j} & =-\int_{D}\Psi_{i}(x)p(x)dx\quad\text{for}\quad i\neq j.
\end{cases}
\]
 An easy computation shows that its eigenvalues are all but one equal
to one. The remaining eigenvalue is 
\[
\lambda_{N}=1-\sum_{i=1}^{N}\int_{D}\Psi_{i}(x)p(x)dx
\]
 and we have 
\[
1-\lambda_{N}=\sum_{i=1}^{N}\left(\int_{D}\varphi_{i}(x)p(x)dx-\frac{\frac{1}{M}\sum_{j=1}^{M}\varphi_{i}(Y_{j})}{\frac{1}{M}\sum_{j=1}^{M}\varphi_{i_{0}}(Y_{j})}\int_{D}\varphi_{i_{0}}(x)p(x)dx\right)
\]
 and the inequality 
\begin{eqnarray*}
\left|1-\lambda_{N}\right|\left|\frac{1}{M}\sum_{j=1}^{M}\varphi_{i_{0}}(Y_{j})\right| & \leq & \left|\int_{D}\sum_{i=1}^{N}\varphi_{i}(x)p(x)dx-\frac{1}{M}\sum_{j=1}^{M}\sum_{i=1}^{N}\varphi_{i}(Y_{j})\right|\left|\int_{D}\varphi_{i_{0}}(x)p(x)dx\right|\\
 &  & +\left|\int_{D}\varphi_{i_{0}}(x)p(x)dx)-\frac{1}{M}\sum_{j=1}^{M}\varphi_{i_{0}}(Y_{j})\right|\left|\int_{D}\sum_{i=1}^{N}\varphi_{i}(x)p(x)dx\right|
\end{eqnarray*}
 which proves that $\lambda_{N}$ converges to 1 when $M\rightarrow\infty$
and that the condition number is once again asymptotically one. Note
that we compute with this formulation an approximation 
\[
P_{N}u(x)=\sum_{i=1}^{N}u(x_{i})\Psi_{i}(x)
\]
 of the solution with discrete integral equal to zero with respect
to the particle approximation of $p(x).$

\subsection{Application to the Poisson equation }

In this section, we describe the application of the stochastic spectral
formulation on our main example of the Poisson equation with pure
Neumann boundary conditions studied in section \ref{sec:Numerical-results}.
The basis functions will rely on Tchebychef interpolation polynomials
in dimension 2. The big advantage of this test case is that the centering
procedure can be done exactly as we know that the invariant probability
is the uniform law and because the integration domain is a square.
This also enables us to compare the exact centering procedure and
the approximate one where the numerical integration is done by means
of a particle approximation of the invariant measure.

\subsubsection{Basis functions}

Our spectral approximation is based on the standard tensorized interpolation
of the solution on the Tchebychef grid which writes 
\[
\sum_{n=0}^{N}\sum_{m=0}^{N}u_{n,m}l_{n}(x)l_{m}(y)
\]
 where $l_{n}$ is the Lagrange polynomial associated to $z_{n}$
and $u_{n,m}$ is the approximate value of the solution at the point
$(z_{n},z_{m})$ where $z_{n}=\cos(\frac{2n+1}{2N+2}\pi$), $0\leq n\leq N$.
For the sake of simplicity, we choose $N$ even so that point $(0,0)$
belongs to the Tchebychef grid. Indeed, in this case the maximum of
the integrals 
\[
\int_{-1}^{1}\int_{-1}^{1}l_{n}(x)l_{m}(y)dxdy
\]
 is always attained for the Lagrange polynomials $l_{\frac{N}{2}}(x)$
and $l_{\frac{N}{2}}(y)$ corresponding to this point. The function
$l_{\frac{N}{2}}(x)l_{\frac{N}{2}}(y)$ is removed from the basis
functions and the $(N+1)^{2}-1$ centered basis functions $\Psi_{n,m}(x,y)$
now write 
\[
\Psi_{n,m}(x,y)=l_{n}(x)l_{m}(y)-\frac{\int_{-1}^{1}\int_{-1}^{1}l_{n}(x)l_{m}(y)dxdy}{\int_{-1}^{1}\int_{-1}^{1}l_{\frac{N}{2}}(x)l_{\frac{N}{2}}(y)dxdy}l_{\frac{N}{2}}(x)l_{\frac{N}{2}}(y)
\]
 for $(n,m)\neq(\frac{N}{2},\frac{N}{2}).$ In the case of the approximate
centering, the integrals above are replaced by their particle approximations.

\subsubsection{Numerical results: exact centering}

In table \ref{tab:Global-error-of}, we present our results based
on the Euler scheme approximation with two time discretization parameters
$\delta_{1}=0.01,\delta_{2}=0.001,$ a regularization parameter $\xi=0.001$
and with two different numbers of simulations $M_{1}=1000$ and $M_{2}=5000$.
The trajectories are stopped at final time $T_{0}=10.$ For these
four sets of parameters, we compute the maximum absolute error over
the grid points 
\[
err_{1}(N)=\max_{0\leq i,j\leq N}\left|u(z_{i},z_{j})-u_{i,j}\right|
\]
 and the condition number $\kappa(N)$ of the spectral matrix for
$N_{1}=2$ and $N_{2}=4$ and $\alpha=\frac{1}{3}.$

\begin{table}
\begin{centering}
\begin{tabular}{|c|c|c|c|c|}
\hline 
 & $err_{1}(N_{1})$  & $\kappa(N_{1})$  & $err_{1}(N_{2})$  & $\kappa(N_{2})$\tabularnewline
\hline 
$(\delta_{1},M_{1})$  & 6.1E-3  & 11.4  & 3.6E-5  & 856\tabularnewline
\hline 
$(\delta_{1},M_{2})$  & 3.4E-3  & 2.9  & 1.2E-5  & 121\tabularnewline
\hline 
$(\delta_{2},M_{1})$  & 4.2E-3  & 5.1  & 2.3E-5  & 114\tabularnewline
\hline 
$(\delta_{2},M_{2})$  & 1.1E-3  & 2.3  & 3.1E-6  & 17\tabularnewline
\hline 
\end{tabular}
\par\end{centering}

\caption{Global error $err_{1}(N)=\max_{0\leq i,j\leq N}\left|u(z_{i},z_{j})-u_{i,j}\right|$of
the approximation of the solution of the pure Neumann problem with
Euler scheme with time steps $\delta_{1}=0.01,\delta_{2}=0.001,$
and a regularization parameter $\xi=0.001$ for for $N_{1}=2$ and
$N_{2}=4$ and $\alpha=\frac{1}{3}.$ The value of condition number
of the spectral matrix is $\kappa(N)$. \label{tab:Global-error-of}}
\end{table}

We can observe that the condition number $\kappa$ is decreasing as $M$
increases and $\delta$ decreases. The system is very well conditioned
especially for the parameters $(\delta_{2},M_{2}).$ As for usual
standard spectral methods, the error is small and decreases with $\kappa$
and when $N$ increases.

\subsubsection{Numerical results: approximate centering}

The spectral method now requires the approximation of the invariant
measure which should be done with $Q$ simulations and a time step
$\delta.$ The invariant measure is the uniform law in the square
$[-1,1]^{2}$. To study its impact on the spectral matrix, we have
chosen in table \ref{tab:Global-approximation-of-1} to make its approximation
using simply $Q$ samples $(X_{l},Y_{l})$ of this uniform law. We
use two samples of different sizes $Q_{1}=100$ and $Q_{2}=10000$
and introduce a new error criterion 
\[
err_{2}(N)=\max_{0\leq i,j\leq N}\left|u(z_{i},z_{j})-u_{i,j}+\frac{1}{4Q}\sum_{l=1}^{Q}u(X_{l},Y_{l})\right|
\]
 as we have proven that we approximate the solution with discrete
integral equal to zero with respect to the particle approximation
of $p(x).$ We nevertheless keep also the previous error criterion
to study the impact of the particle approximation on the bias.

\begin{table}
\begin{centering}
\begin{tabular}{|c|c|c|c|c|c|c|}
\hline 
 & $err_{1}(N_{1})$  & $err_{2}(N_{1})$  & $\kappa(N_{1})$  & $err_{1}(N_{2})$  & $err_{2}(N_{2})$  & $\kappa(N_{2})$\tabularnewline
\hline 
$(\delta_{1},M_{1},Q_{1})$  & 1.9E-2  & 2.4E-3  & 10.2  & 4.60E-2  & 2.7E-5  & 550\tabularnewline
\hline 
$(\delta_{1},M_{1},Q_{2})$  & 5.8E-3  & 5.1E-3  & 26.3  & 1.5E-3  & 2.3E-4  & 1187\tabularnewline
\hline 
$(\delta_{1},M_{2},Q_{1})$  & 1.7E-2  & 1.6E-3  & 5.2  & 1.4E-2  & 6.3E-5  & 1130\tabularnewline
\hline 
$(\delta_{1},M_{2},Q_{2})$  & 5.1E-3  & 6.1E-4  & 2.9  & 3.1E-3  & 1.0E-5  & 644\tabularnewline
\hline 
$(\delta_{2},M_{1},Q_{1})$  & 1.6E-2  & 3.6E-3  & 4.1  & 1.5E-2  & 5.4E-5  & 309\tabularnewline
\hline 
$(\delta_{2},M_{1},Q_{2})$  & 3.9E-3  & 1.4E-3  & 3.0  & 7.4E-4  & 1.4E-5  & 252\tabularnewline
\hline 
$(\delta_{2},M_{2},Q_{1})$  & 6.3E-2  & 9.4E-4  & 2.8  & 3.0E-2  & 3.1E-6  & 22\tabularnewline
\hline 
$(\delta_{2},M_{2},Q_{2})$  & 1.5E-3  & 1.2E-3  & 1.6  & 3.3E-3  & 4.2E-6  & 16\tabularnewline
\hline 
\end{tabular}
\par\end{centering}

\caption{Global approximation of the solution of the pure Neumann problem with
an Euler scheme. The grid parameters are respectively $N_{1}=2$ and
$N_{2}=4$. The time steps are $\delta_{1}=0.01,\delta_{2}=0.001$.
The invariant measure is approximated with samples of sizes $Q_{1}=100$
and $Q_{2}=10000$. \label{tab:Global-approximation-of-1}}
\end{table}

The key observation is that we check that the solution effectively
computed is the one with discrete integral equal to zero with respect
to the particle approximation of $p(x).$ Both accuracy and condition
number have the same behaviour with respect to the parameters $\delta$
and $M$ than with the exact centering. The condition number and $err_{1}$
decrease when the particle approximation is done with more points.

\subsection{Application to a convection-diffusion problem}

In this section, we study a more general convection diffusion equation
with an additional drift coefficient depending assymmetrically on
the spatial position. We keep nevertheless a square domain in order
not to mix the approximation problems in a general bounded domain
and the other difficulties met by our method. We perform our algorithm
on the operator

\[
L=\frac{1}{2}\Delta+x\beta_{x}\frac{\partial}{\partial x}+y\beta_{y}\frac{\partial}{\partial y}
\]
with a source term $f(x,y)=\alpha(\alpha+x\beta_{x}+y\beta_{y})\exp(\alpha(x+y))$
and a boundary term $g(x,y)=\pm\frac{\alpha}{2}\exp(\alpha(x+y))$
on $\partial D$ such that the solution up to an additive constant
is still $\exp(\alpha(x+y)).$ For this model, we do not know the
exact invariant measure of the reflected diffusion process associated
to the operator $L$. To compute its particle approximation, we run
a single path starting at the center of the domain with a time step
$\delta$ and ending at time $T=Q\delta.$ The $Q$ points of this
path give an approximation of the invariant measure. We present in
table \ref{tab:Global-approximation-of} our results for the parameters
$\beta_{x}=0.2,\beta_{y}=0.1$ and $\alpha=0.3$. We observe that
the accuracy of the algorithm is still very high, the condition number
of the system still very low even if we approximate the invariant
measure by running a path. This confirms that our approach is quite
general.

\begin{table}
\begin{centering}
\begin{tabular}{|c|c|c|c|c|}
\hline 
 & $err_{2}(N_{1})$  & $\kappa(N_{1})$  & $err_{2}(N_{2})$  & $\kappa(N_{2})$\tabularnewline
\hline 
$(\delta_{1},M_{1},Q_{1})$  & 4.7E-3  & 16.1  & 2.E-5 & 2648\tabularnewline
\hline 
$(\delta_{1},M_{1},Q_{2})$  & 2.6E-3  & 26.4 & 1.1E-5  & 649\tabularnewline
\hline 
$(\delta_{1},M_{2},Q_{1})$  & 1.1E-2  & 74.2  & 1.0E-5  & 1254\tabularnewline
\hline 
$(\delta_{1},M_{2},Q_{2})$  & 3.1E-3  & 12.0  & 4.3E-4  & 22437\tabularnewline
\hline 
$(\delta_{2},M_{1},Q_{2})$  & 2.7E-3  & 5.9  & 2.1E-5  & 1200\tabularnewline
\hline 
$(\delta_{2},M_{2},Q_{2})$  & 8.1E-4  & 4.0  & 2.0E-6 & 39.8\tabularnewline
\hline 
\end{tabular}
\par\end{centering}

\caption{Global approximation of the solution of the pure Neumann problem with
an Euler scheme. The grid parameters are respectively $N_{1}=2$ and
$N_{2}=4$. The time steps are $\delta_{1}=0.01,\delta_{2}=0.001$.
The invariant measure is approximated with samples of sizes $Q_{2}=10000$
and $Q_{1}=100$ only for $\delta_{1}$. The final time $Q_{1}\delta_{2}=0.1$
is not large enough to reasonnably approximate the invariant measure.
\label{tab:Global-approximation-of}}
\end{table}

\section{Conclusion\label{sec:Conclusion}}

To compute Monte Carlo approximations of the solution of the Neumann
problem for elliptic equations, we had to overcome several difficulties.
First, we have characterized the solution of the Feynman-Kac representation
introduced in \cite{bencherif_pardoux} as the one with zero mean
value with respect to the invariant measure of its associated stochastic
process. Then, we have proven that the variance increases mainly linearly
as a function of the time $T_{0}$ to stop the trajectories.

We have introduced some new schemes to deal with the inhomogeneous
Neumann boundary conditions. They were tested first on a pointwise
approximations of mixed Dirichlet-Neumann problem where they show
a good efficiency. The pointwise resolution of the pure Neumann problem
was a lot harder. Indeed, we had to choose $T_{0}$ not too large
because of the increase of the variance but also not too small because
of the bias. We had also to understand that the solution computed
numerically depends on the parameters of the numerical schemes and
is not equal to the one with zero mean value with respect to the invariant
measure. Taken all these difficulties into account, we have been able
nevertheless to reach a reasonable accuracy on the approximate solutions.

Concerning the global spectral approximation, we had to pay attention
to the zero mean value property of the solution to chose our approximation
basis. This has been achieved using exact or approximate centering
procedures very similar to the usual ones used in finite element methods.
In both cases, the condition number of the spectral matrix was proven
to be asymptotically one. The numerical experiments show that the
stochastic spectral method is both very accurate and well-conditioned.

The pointwise approximations of the pure Neumann problem are not completely
satisfactory because the solutions obtained depend on the parameters
of the numerical scheme. The choice of the time $T_{0}$ to stop the
trajectories is not straightforward for a general diffusion in a complex
domain. For the WOS method, it requires furthermore to keep in memory
discretisation of trajectories which is both costly and adds an error
not so easy to quantify. It could be interesting to add a penalization
term either in the source term or via Robin boundary conditions to
at least get rid of some of these drawbacks. 

\section*{Acknowledgments}
We would like to thank the Anonymous Referee for his/her 
constructive comments on our work and for drawing our attention to  some references we were not aware of.

\end{document}